\documentclass[12pt]{article}
\usepackage{amsmath, amsthm, amssymb,latexsym,enumerate}
\usepackage{graphicx}
\usepackage{hyperref}
\usepackage{xcolor}
\usepackage{todonotes}
\usepackage{cite}
\usepackage{dsfont}

\newcommand{\vp}{\varphi}
\newcommand{\uc}{\vp_o}

\newtheorem{theorem}{Theorem}[section]
\newtheorem{lemma}[theorem]{Lemma}
\newtheorem{corollary}[theorem]{Corollary}
\newtheorem{claim}[theorem]{Claim}

\newtheorem{conjecture}[theorem]{Conjecture}
\newtheorem{question}[theorem]{Question}

\numberwithin{equation}{section}

\DeclareMathOperator{\mad}{mad}

\title{Odd coloring of sparse graphs and planar graphs}

\author{Eun-Kyung Cho\thanks{
Department of Mathematics, Hankuk University of Foreign Studies, Yongin-si, Gyeonggi-do, Republic of Korea.
 \texttt{ekcho2020@gmail.com}
}
\and Ilkyoo Choi\thanks{
Department of Mathematics, Hankuk University of Foreign Studies, Yongin-si, Gyeonggi-do, Republic of Korea.
\texttt{ilkyoo@hufs.ac.kr}
}
\and Hyemin Kwon\thanks{
Department of Mathematics, Ajou University, Suwon-si, Gyeonggi-do, Republic of Korea.
\texttt{khmin1121@ajou.ac.kr}
}
\and Boram Park\thanks{
Department of Mathematics, Ajou University, Suwon-si, Gyeonggi-do, Republic of Korea.
\texttt{borampark@ajou.ac.kr}
}}
\date{}

\begin{document}
 
\maketitle

\begin{abstract} 
An {\it odd $c$-coloring} of a graph is a proper $c$-coloring such that each non-isolated vertex has a color appearing an odd number of times on its neighborhood.
This concept was introduced very recently by Petru\v sevski and \v Skrekovski and has attracted considerable attention.  
Cranston investigated odd colorings of graphs with bounded maximum average degree, and conjectured that every graph $G$ with $\mad(G)\leq \frac{4c-4}{c+1}$ has an odd $c$-coloring for $c\geq 4$, and proved the conjecture for $c\in\{5, 6\}$. 
In particular, planar graphs with girth at least $7$ and $6$ have an odd $5$-coloring and an odd $6$-coloring, respectively. 

We completely resolve Cranston's conjecture.
For $c\geq 7$, we show that the conjecture is true, in a stronger form that was implicitly suggested by Cranston, but for $c=4$, we construct counterexamples, which all contain $5$-cycles. 
On the other hand, we show that a graph $G$ with $\mad(G)<\frac{22}{9}$ and no induced $5$-cycles has an odd $4$-coloring. 
This implies that a planar graph with girth at least 11 has an odd $4$-coloring.
We also prove that a planar graph with girth at least 5 has an odd $6$-coloring.
\end{abstract}

\section{Introduction}

All graphs in this paper are finite and simple, which means no loops and no parallel edges. 
For a graph $G$, let $V(G)$ and $E(G)$ denote its vertex set and edge set, respectively.
An {\it odd $c$-coloring} of a graph is a proper $c$-coloring  with the additional constraint that  each  vertex of positive degree has a color appearing an odd number of times among its neighbors. 
A graph is {\it odd $c$-colorable} if it has an odd $c$-coloring.
The {\it odd chromatic number} of a graph $G$, denoted $\chi_o(G)$, is the minimum $c$ such that $G$ has an odd $c$-coloring.  

This concept was introduced very recently by Petru\v sevski and \v Skrekovski~\cite{petrusevski2021colorings}, who showed that planar graphs are odd 9-colorable, and conjectured that planar graphs are odd 5-colorable.
Note that a 5-cycle is a planar graph whose odd chromatic number is exactly 5. 
A study on various aspects of the odd chromatic number was carried out by Caro, Petru\v sevski, and \v Skrekovski~\cite{caro2022remarks}.
Building upon some results in~\cite{caro2022remarks}, Petr and Portier~\cite{petr2022odd} proved that planar graphs are odd 8-colorable. 
A strengthening of the aforementioned result was shown by Fabrici et al.~\cite{fabrici2022proper}, where they also compiled results on planar graphs regarding similar coloring parameters, which may be of independent interest~\cite{2020BhKa,2018AbAlDeFeGoHeKeSc}.
Odd colorings of $1$-planar graphs and $k$-planar graphs were also studied, see~\cite{cranston2022note,DuMoOd,Hickingbotham}.

Cranston~\cite{cranston2022odd} investigated odd colorings of sparse graphs, measured in terms of the maximum average degree. 
Given a graph $G$, the {\it maximum average degree} of $G$, denote $\mad(G)$, is the maximum of $\frac{2|E(H)|}{|V(H)|}$ over all non-empty subgraphs $H$ of $G$. 

For $c\in\{1,2\}$, it is not hard to see that maximum average degree at most $c-1$ guarantees odd $c$-colorability, and the bound on maximum average is tight. 
Also, for $c\in\{3,4\}$, maximum average degree less than $2$ guarantees odd $c$-colorability, and the bound on the maximum average degree cannot be changed to include equality because of the $5$-cycle. 
See also~\cite{cranston2022odd}. 
For $c\geq 4$, Cranston~\cite{cranston2022odd} made the following conjecture:

\begin{conjecture}[\cite{cranston2022odd}]\label{conj:cranston}
For $c\geq 4$, if $G$ is a graph with $\mad(G) < \frac{4c}{c+2}$, then $\chi_o(G)\leq c$.
\end{conjecture}

Let $K^*_{c+1}$ be the graph obtained by subdividing every edge of the complete graph on $c+1$ vertices exactly once. 
Note that $\mad(K^*_{c+1})=\frac{4c}{c+2}$ and $\chi_o(K^*_{c+1})=c+1$ if $c\neq 1$.

Cranston proved a weakening of the conjecture by showing that a graph $G$ with $\mad(G)\leq\frac{4c-4}{c+1}$ is odd $(c+3)$-colorable. 
He also proved the conjecture for $c\in\{5,6\}$ in a stronger form.

We completely resolve Conjecture~\ref{conj:cranston}.
For $c\geq7$, we confirm Conjecture~\ref{conj:cranston} in the affirmative, but for $c=4$, we construct counterexamples. 
Actually, we prove the below stronger statement, which was implicitly suggested in~\cite{cranston2022odd}:

\begin{theorem}\label{thm:mad}
For $c\geq 7$, if $G$ is a graph with $\mad(G)\leq\mad(K^*_{c+1})=\frac{4c}{c+2}$, then $\chi_o(G)\leq c$, unless $G$ contains $K^*_{c+1}$ as a subgraph.
\end{theorem}

We remark that Cranston confirmed Conjecture~\ref{conj:cranston} for $c\in\{5,6\}$ in the form of Theorem~\ref{thm:mad}.

For $c=4$, Conjecture~\ref{conj:cranston} is false as $\mad(C_5)=2<\frac{8}{3}$ yet $\chi_o(C_5)=5>4$. 
Moreover, there are actually infinitely many counterexamples:
Let $H_k$ be the graph obtained by choosing one vertex $x_i$ each from $k$ disjoint $5$-cycles, then identifying all $x_i$ into a single vertex $v$.
If there are only four colors, then the two neighbors of $v$ from the same $5$-cycle must be of the same color.
Hence, every color appears an even number of times among the neighbors of $v$, so $H_k$ is not odd 4-colorable.
Since $\mad(H_k)=\frac{10k}{4k+1}<\frac{8}{3}$, each $H_k$ is a  counterexample for all $k \ge 1$.
Note that $H_1$ is the $5$-cycle.
Furthermore, Caro, Petru\v sevski, and \v Skrekovski~\cite{caro2022pcfremarks} showed that every non-trivial connected graph where every block of is isomorphic to $C_5$ has odd chromatic number $5$.

Since Conjecture~\ref{conj:cranston} is false for $c=4$ and all counterexamples we found contain an induced $5$-cycle, we investigate the class of graphs with bounded maximum average degree and no induced 5-cycle, and obtain the following result:

\begin{theorem}\label{thm:mad22/9}
If $G$ is a graph with $\mad(G)<\frac{22}{9}$ and no induced 5-cycle, then $\chi_o(G)\leq 4$. 
\end{theorem}

Note that the threshold on the maximum average degree is higher than the one in Conjecture~\ref{conj:cranston}. 

Results regarding bounded maximum average degree have natural corollaries to planar graphs with girth restrictions. 
Namely, the class of graphs $G$ with $\mad(G)<\frac{2g}{g-2}$ includes all planar graphs with girth at least $g$. 
Hence, by Cranston's results, Theorem~\ref{thm:mad}, and Theorem~\ref{thm:mad22/9}, we have the following corollaries:

\begin{corollary}\label{cor:planar5}
For $c\geq 5$, every planar graph with girth at least $\left\lceil \frac{4c}{c-2} \right\rceil$ is odd $c$-colorable. 
\end{corollary}

\begin{corollary}
Every planar graph with girth at least $11$ is odd $4$-colorable. 
\end{corollary}

Corollary~\ref{cor:planar5} states that a planar graph with girth at least $5$ is odd $10$-colorable, and a planar graph with girth at least $6$ is odd $6$-colorable. 
We improve upon both statements by proving the following result: 

\begin{theorem}\label{thm:6colors}
Every planar graph $G$ with girth at least $5$ is odd $6$-colorable. 
\end{theorem}

The paper is organized as  follows.
In Section~\ref{sec:keylem}, we first prove a key lemma that generates reducible configurations for Theorem~\ref{thm:mad} and Theorem~\ref{thm:6colors}.
The discharging procedures used in the proofs for Theorem~\ref{thm:mad} and Theorem~\ref{thm:6colors} appear in Section~\ref{sec:mad} and Section~\ref{sec:6colors}, respectively. 
In Section~\ref{sec:mad22/9}, we provide both the reducible configurations and the discharging procedure for Theorem~\ref{thm:mad22/9}. 
We conclude the paper with possible future research directions in Section~\ref{sec:future}.

\section{A key lemma and preliminaries}\label{sec:keylem}

We first lay out some frequently used notation and terminology.
Given a graph $G$ and $S \subseteq V(G)$, let $G-S$ denote the subgraph of $G$ induced by $V(G) \setminus S$.
For a vertex $v$, let $N_G(v)$ denote the neighborhood of $v$  and $d_G(v)$ denote the degree of $v$.
A $d$-vertex (resp. $d^+$-vertex and $d^-$-vertex) is a vertex of degree exactly $d$ (resp. at least $d$ and at most $d$). 
A $d$-neighbor of a vertex $v$ is a $d$-vertex that is a neighbor of $v$.
Let $N_{G,d}(v)$ be the set of all $d$-neighbors of $v$, and 
$n_{G,d}(v)=|N_{G,d}(v)|$.
We also define $d^+$-neighbor, $d^-$-neighbor, $N_{G, d^+}(v)$, $N_{G, d^-}(v)$, $n_{G, d^+}(v)$, and $n_{G, d^-}(v)$ analogously. 
If there is no confusion, then we drop $G$ in our notation, as in  $d(v)$, $N(v)$, $N_d(v)$, $n_{d^+}(v)$, and so on. 
A {\it $k$-thread} means a path on (at least) $k$ $2$-vertices. 
 
An {\it odd color} of a vertex $v$ is a color that appears an odd number of times on $N_G(v)$.  
Given an odd coloring $\varphi$ of a graph $G$ and a vertex $v$, let $\varphi_o(v)$ denote an odd color of $v$; if $v$ has many odd colors, then choose an arbitrary one. 
Oftentimes we will be extending a partial coloring $\varphi$ of $G$ (to the whole graph $G$), and we will abuse notation and use $\varphi_o(v)$ to denote an odd color of the extended coloring as well. 
Note that in a proper coloring of $G$, a vertex of odd degree will always have an odd color.
For two vertices $x$ and $y$, when we say ``color $x$ with a color that is not $\varphi(y)$ (or $\varphi_o(y)$)'', we exclude the color $\varphi(y)$ (or $\varphi_o(y)$) only when $\varphi(y)$ (or $\varphi_o(y)$) is defined. 

In addition, for technical reasons, we allow an odd $c$-coloring of the null graph, which is a graph whose vertex set is empty. 

The following is a key lemma that is useful in multiple situations.  
A $3^+$-vertex is {\it easy} if it has either odd degree or a $2$-neighbor.
For each vertex $v$, let $N_e(v)$ denote the set of easy neighbors of $v$ and let $n_e(v)=|N_e(v)|$.

\begin{lemma}\label{keylem}
For $c\geq 5$, let $G$ be a graph with no odd $c$-coloring, but every proper subgraph of $G$ has an odd $c$-coloring.  Then $G$ has no $1^-$-vertex, and 
$2d(v)\geq n_2(v)+n_e(v)+c$ for an easy vertex $v$.
\end{lemma}
\begin{proof} 
Since $c\geq 3$, $G$ has no $1^-$-vertex $z$ because an odd $c$-coloring of $G-\{z\}$ can be easily extended to an odd $c$-coloring of $G$, which is a contradiction. Thus $G$ has no $1^-$-vertex. 
Let $X=N(N_2(v))\setminus\{v\}$.

Suppose to the contrary that $2d(v)\leq n_2(v)+n_e(v)+c-1$.
 By the minimality of $G$, there is an odd $c$-coloring $\varphi$ of $G-( N_2(v) \cup \{v\})$.
Use a color not in
$C= \varphi( X\cup N_{3^+}(v)) \cup \varphi_o(N_{3^+}(v)\setminus N_e(v))$ on $v$; this is possible because
\[|C|\leq
|X \cup N_{3^+}(v)| + |N_{3^+}(v)\setminus N_e(v)|\le 
d(v)+(d(v)-n_2(v)-n_e(v)) = 2d(v)-n_2(v)-n_e(v) \leq c-1.\] 
Now, start coloring the vertices in $N_2(v)$ one by one. 
For each $2$-vertex $u$ in $N_2(v)$,
$|\vp(N(u))\cup\uc(N(u))|\leq 4\leq c-1$,
so we can use a color not in $\vp(N(u))\cup\uc(N(u))$ on $u$.
At this point, we have a proper coloring of $G$, and since either $N_2(v)\neq \emptyset$ or $d(v)$ is odd, every vertex of $G$ has an odd color, except maybe the vertices in $N_e(v)$.
If a vertex $u\in N_e(v)$ does not have an odd color, then $d(u)$ is even, so $u$ must have a $2$-neighbor. 
For each such $u$, recolor a $2$-neighbor $x$ of $u$ with a color not in $\vp(N(x))\cup\uc(N(x))$ to guarantee that $u$ has an odd color. 
Note that $|\vp(N(x))\cup\uc(N(x))|\leq 4\leq c-1$.
This is an odd $c$-coloring of $G$, which is a contradiction. 
\end{proof}

\section
{Graphs with bounded maximum average degree}
\label{sec:mad}

Fix $c\geq 7$, and let $G$ be a counterexample to Theorem~\ref{thm:mad} with the minimum number of vertices. 
So, $G$ is a graph with $\mad(G)\leq \frac{4c}{c+2}$ not containing $K^*_{c+1}$ as a subgraph, where $G$ has no odd $c$-coloring. 
Let $\mu(v):=d(v)$ be the initial charge of a vertex $v$. 
Let $\mu^*(v)$ be the final charge of $v$ after the following discharging rule:
\begin{enumerate}[(R1)]
    \item\label{madrule:1} 
    Each $4^+$-vertex sends charge $\frac{c-2}{c+2}$ to each of its $2$-neighbors.
\end{enumerate}

Note that the total charge is preserved during the discharging process.  Consider a vertex $v$.  
We will show that $\mu^*(v) \geq \frac{4c}{c+2}$.  
Since $c\ge 7$, Lemma~\ref{keylem} implies that there is neither a $1$-vertex nor a $3$-vertex.

If $v$ is a 2-vertex, then it receives charge $\frac{c-2}{c+2}$ from each of its neighbors by (R\ref{madrule:1}) since $v$ has only $4^+$-neighbors.
Note that $v$ has no $2$-neighbors by Lemma~\ref{keylem}.
Hence,  $\mu^*(v)=2+\frac{c-2}{c+2}\cdot 2=\frac{4c}{c+2}$.
If $v$ is an easy $4^+$-vertex, then $n_2(v)\le \min\{d(v), 2d(v)-c\}$ by Lemma~\ref{keylem}, so by~(R\ref{madrule:1}), 
$\mu^*(v)\geq 
 d(v)-\frac{c-2}{c+2}\cdot \min\{d(v), 2d(v)-c\}$. 
Then \[\mu^*(v)\ge \begin{cases}
\begin{aligned}
 & d(v)-\frac{c-2}{c+2}\cdot d(v)  
 = \frac{4d(v)}{c+2} 
 \ge\frac{4c}{c+2},
 &\mbox{ if $d(v)\geq c$.} \\
 & d(v)-\frac{c-2}{c+2}\cdot(2d(v)-c) 
 = \frac{-(c-6)d(v)+c(c-2)}{c+2}
 >\frac{4c}{c+2},
 &\mbox{ if $d(v)< c$.}
 \end{aligned}
\end{cases}\] 
If $v$ is a non-easy $4^+$-vertex, then it has no $2$-neighbors, so $v$ is not involved in the discharging rule. 
Hence, $\mu^*(v)=d(v) \ge 4 >\frac{4c}{c+2}$.

Since $\mad(G) \le \frac{4c}{c+2}$, the initial charge sum is at most $\frac{4c}{c+2}|V(G)|$. 
By the above argument, the final charge sum is at least $\frac{4c}{c+2}|V(G)|$.
Henceforth, $\mu^*(v)=\frac{4c}{c+2}$ for every vertex $v$.  

From the discharging procedure, it follows that $G$ has only $2$-vertices and $c$-vertices, and every $c$-vertex has only $2$-neighbors and vice versa. 
Thus, $G$ is obtained from a connected $c$-regular graph $G_0$ by subdividing every edge of $G_0$. 
Since $G_0$ is not $K_{c+1}$, we know that $G_0$ has a proper $c$-coloring $\varphi$ by Brooks' Theorem~\cite{1941Brooks}.
Consider $\varphi$ as a coloring of the $c$-vertices of $G$, and then  color the 2-vertices of $G$ one by one as follows. 
For each $2$-vertex $v$, color $v$ by a color not in $\vp(N(v))\cup\uc(N(v))$.
Now, $\varphi$ is an odd $c$-coloring of $G$, which is a contradiction.

\section{Planar graphs with girth at least $5$ are odd $6$-colorable}
\label{sec:6colors}

Let $G$ be a counterexample to Theorem~\ref{thm:6colors} with the minimum number of vertices, and fix a plane embedding of $G$. 
So, $G$ is a plane graph with girth at least $5$ that has no odd $6$-coloring.

\begin{figure}[h!]\centering
\includegraphics[width=15cm]{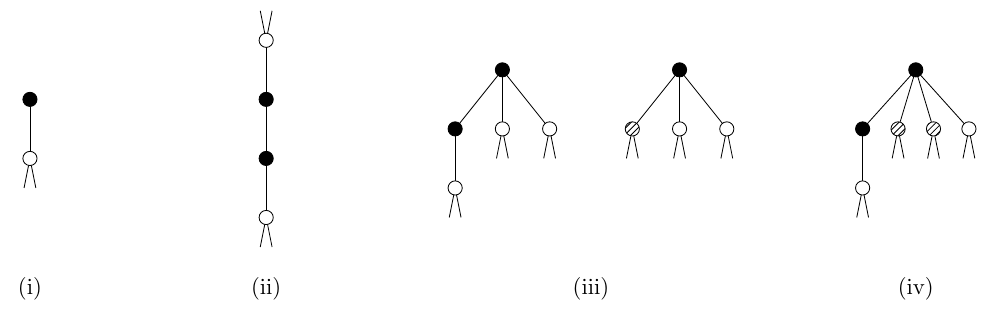}
\caption{An illustration for Claim~\ref{6clm:struc}.
Hatched vertices are easy vertices and black vertices have all its incident edges drawn.
% White vertices may have additional incident edges.
}
\label{fig1}
\end{figure}

\begin{claim}\label{6clm:struc}
The following do not appear in $G$:
\begin{enumerate}[\rm (i)]
    \item a $1$-vertex,
    \item a $2$-thread,   
    \item a $3$-vertex with either  a $2$-neighbor or an easy neighbor, and
    \item an easy $4$-vertex with two easy neighbors.
\end{enumerate}
\end{claim}
\begin{proof}
Each of (i), (ii), (iii), (iv) follows from Lemma~\ref{keylem}. See Figure~\ref{fig1} for an illustration.
\end{proof}

A face is {\it bad} if it is a $5$-face incident with a $2$-vertex, and a face is {\it good} if it is a $5$-face that is not bad.
 
For a face $f$, let $d(f)$ denote the degree of $f$, which is the length of a boundary walk of $f$.
For a vertex or a face $z$, define the initial charge $\mu(z)$ of $z$ as the following: for each vertex $v$ let $\mu(v)=d(v)-6$, and for each face $f$ let $\mu(f)=2d(f)-6$.
Let $\mu^*(z)$ be the final charge of $z$ after the following discharging rules:

\begin{enumerate}[(R1)]
    \item\label{6rule:2vx} Each face sends charge $2$ to each incident $2$-vertex.
    \item\label{6rule:3vx} Each face sends charge $1$ to each incident $3$-vertex.
    \item\label{6rule:bad5} Each bad face $f$ incident with exactly one $2$-vertex sends charge $\frac{1}{2}$ to each incident $4^+$-vertex that has a $4^+$-neighbor on $f$. 
    \item\label{6rule:good5} Each good face sends charge $1$ to each incident easy $4$-vertex, and sends charge $\frac{1}{2}$ to each incident $4^+$-vertex that is not an easy $4$-vertex.      
    \item\label{6rule:6face} For a $4^+$-vertex $v$ on a $6^+$-face $f$,
    if $xvz$ are three consecutive vertices on a boundary walk of $f$, then $f$ sends charge $\frac{1}{2} \cdot |\{x,z\} \setminus N_2(v)|$ to $v$.  
\end{enumerate}

Note that the total charge is preserved during the discharging process.

\begin{claim}
Each face has non-negative final charge.
\end{claim}
\begin{proof}
Suppose $f$ is a $6^+$-face.
We first send charge $1$ to each incident edge, which is possible since $\mu^*(f)=2d(f)-6-d(f)\cdot1\geq 0$, and $d(f)\geq 6$.
If an edge $e$ is incident with a $2$-vertex $x$, then $e$ sends charge $1$ to $x$, 
and if $e$ is not incident with a $2$-vertex, then $e$ sends charge $\frac{1}{2}$ to each vertex incident with $e$. 
Consider a vertex $v$. 
If $v$ is a $2$-vertex, then it will receive charge $2$ from $f$ since both edges incident with $v$ will send charge $1$ to $v$. 
If $v$ is a $3$-vertex, then it will receive charge $1$ from $f$ since $v$ is not adjacent to a $2$-vertex by Claim~\ref{6clm:struc} (iii).
If $v$ is a $4^+$-vertex such that $xvz$ are three consecutive vertices on a boundary walk of $f$, then $v$ will receive charge $\frac{1}{2} \cdot |\{x,z\} \setminus N_2(v)|$ from $f$.

Suppose $f$ is a bad face. 
Note that $f$ is incident with at most two $3^-$-vertices, which are not adjacent to each other by Claim~\ref{6clm:struc} (ii) and (iii).
If $f$ is incident with exactly two $2$-vertices, then it sends charge $2$ to each of its incident $2$-vertices by (R\ref{6rule:2vx}), so $\mu^*(f)=4-2\cdot2=0$.
Suppose that $f$ has exactly one $2$-vertex. If $f$ is incident with a $3$-vertex, then it sends charge $2$ to its incident $2$-vertex, charge $1$ to its incident $3$-vertex, and charge $\frac{1}{2}$ to two of its incident $4^+$-vertices by (R\ref{6rule:2vx}), (R\ref{6rule:3vx}), and (R\ref{6rule:bad5}).
Hence, $\mu^*(f)= 4-2-1-2\cdot\frac{1}{2}=0$.
If $f$ is not incident with a $3$-vertex, then it sends charge $2$ to its incident $2$-vertex and charge $\frac{1}{2}$ to each of its incident $4^+$-vertices by (R\ref{6rule:2vx}) and (R\ref{6rule:bad5}).
Hence, $\mu^*(f)=4-2-4\cdot\frac{1}{2}=0$.

Suppose, $f$ is a good face, so it is not incident with a $2$-vertex. 
By Claim~\ref{6clm:struc}~(iii) and (iv), $f$ is incident with at most three vertices that are either a $3$-vertex or an easy $4$-vertex.
Thus, $f$ sends charge $1$ to each of its incident $3$-vertices and easy $4$-vertices and charge $\frac{1}{2}$ to each of the other incident $4^+$-vertices by (R\ref{6rule:3vx}) and (R\ref{6rule:good5}).
Hence, $\mu^*(f) \ge 4-3\cdot 1 - 2\cdot \frac{1}{2} =0$.
\end{proof}

\begin{claim}
Each vertex has non-negative final charge.
\end{claim}
\begin{proof}
Consider a vertex $v$.
Note that $v$ is not a $1$-vertex by Claim~\ref{6clm:struc}~(i).
If $v$ is a $2$-vertex, then it always receives charge $2$ twice from its incident faces by (R\ref{6rule:2vx}).
Hence, $\mu^*(v)=-4+2\cdot2=0$.
If $v$ is a $3$-vertex, then it always receives charge $1$ three times from  its incident faces by (R\ref{6rule:3vx}).
Hence, $\mu^*(v)=-3+3\cdot1=0$.

Note that by the rules, a face sends charge either $0$ or at least $\frac{1}{2}$ to its incident $4^+$-vertex.
Moreover, for a $4^+$-vertex $v$ on a face $f$ where $xvz$ are three consecutive vertices on a boundary walk of $f$, if $f$ does not send charge to $v$ by the rules, then one of $x,z$ is a $2$-vertex and the other is either a $3^-$-vertex or an easy vertex.  

Suppose $v$ is a $4$-vertex with neighbors $u_1,u_2,u_3, u_4$ in  cyclic order around $v$. 
If all faces incident with $v$ send charge to $v$, then by the rules, $\mu^*(v) \ge -2+4\cdot\frac{1}{2}=0$.
Suppose that $v$ has an incident face $f$ that does not send charge to $v$, and assume $u_1vu_2$ are three consecutive vertices on a boundary walk of $f$, so $n_2(v)+n_e(v)\geq 2$.
Lemma~\ref{keylem} further implies $n_2(v) + n_e(v) = 2$, so both $u_3$ and $u_4$ are non-easy $4^+$-neighbors. 
Hence, the face incident with $u_3vu_4$ is either a good face or a $6^+$-face, which sends charge $1$ to $v$ by~(R\ref{6rule:good5}) or (R\ref{6rule:6face}), and the faces incident with $u_2vu_3$ and $u_4vu_1$ will each send charge at least $\frac{1}{2}$ to $v$ by (R\ref{6rule:bad5}), (R\ref{6rule:good5}), or (R\ref{6rule:6face}).
Hence, $\mu^*(v) =-2+1+2\cdot\frac{1}{2}=0$.

If $v$ is a $5$-vertex, then $n_2(v)+n_e(v) \le 4$ by Lemma~\ref{keylem}, which further implies that $v$ receives charge at least $\frac{1}{2}$ twice from its incident faces, so by the rules, 
$\mu^*(v) \ge -1 + 2 \cdot \frac{1}{2} = 0$.
 
If $v$ is a $6^+$-vertex, then it sends no charge by the rules, hence, $\mu^*(v)\ge \mu(v)= d(v)-6\geq 0$.
\end{proof}

Therefore, every vertex and every face of $G$ have non-negative final charge, while the initial charge sum is at most $-12$ by the Euler's formula.
This is a contradiction.

\section{Sparse graphs without induced $5$-cycles are odd $4$-colorable}
\label{sec:mad22/9}

Let $G$ be a counterexample to Theorem~\ref{thm:mad22/9} with the minimum number of vertices.
So, $G$ is a graph with $\mad(G) < \frac{22}{9}$ that does not contain an induced $5$-cycle and has no odd $4$-coloring.

Note that we cannot use Lemma~\ref{keylem}, since we have only four colors. Recall that a {\it $k$-thread} is a path on (at least) $k$ $2$-vertices. 

\begin{claim}\label{clm:mad_rc}
The following do not appear in $G$:

\begin{enumerate}[(i)]
    \item a $1$-vertex,
    \item an $\ell$-cycle incident with $(\ell-1)$ $2$-vertices for $\ell\in\{3,4\}$, 
    \item a $4$-thread,
    \item a vertex of odd degree adjacent to a $2$-thread, and
    \item a $3$-vertex with only $2$-neighbors.
\end{enumerate}
\end{claim}
\begin{proof} (i) Suppose to the contrary that $G$ has a $1$-vertex $v$.
Let $u$ be the neighbor of $v$. 
Let $S=\{v\}$, and obtain an odd 4-coloring $\varphi$ of $G-S$. 
Use a color not in $\{\varphi(u), \varphi_o(u)\}$ on $v$ to extend $\varphi$ to all of $G$.

(ii) Suppose to the contrary that $G$ has an $\ell$-cycle $u_1\ldots u_\ell$ of $2$-vertices except $u_1$. 
Let $S = \{u_2, \ldots u_\ell\}$, and obtain an odd $4$-coloring $\varphi$ of $G-S$.
Now, we can obtain an odd 4-coloring of $G$ by greedily coloring $u_2$,  ($u_4$ if it exists,) and $u_3$ in this order. 
Note that when coloring $u_i$, there are at most three colors in  $\{\varphi(u_{i-1}), \varphi_o(u_{i-1}), \varphi(u_{i+1}), \varphi_o(u_{i+1})\}$ (consider the indices modulo $\ell$), so there is a color to use for each $u_i$.

(iii)
Suppose to the contrary that $G$ has a path $u_1v_1v_2v_3v_4u_2$ where $v_1v_2v_3v_4$ is a $4$-thread.
Let $S=\{v_1, v_2, v_3, v_4\}$, and obtain an odd 4-coloring $\varphi$ of $G-S$.
Note that $u_i$'s and $v_i$'s are all distinct since $G$ has no induced $5$-cycle, and by (ii). 
Thus $\varphi(u_1)$, $\varphi(u_2)$, $\varphi_o(u_1)$, and $\varphi_o(u_2)$ are defined, and $v_1$ and $v_4$ have distance at least $3$.

If $\{\varphi(u_1), \varphi_o(u_1)\}$ and $\{\varphi(u_2), \varphi_o(u_2)\}$ are not disjoint,  then there exists a color $c$ that can be used on $v_1$ and $v_4$. 
Now, we can greedily color $v_2$ and $v_3$ to obtain an odd 4-coloring of $G$.

If $\{\varphi(u_1), \varphi_o(u_1)\}$ and $\{\varphi(u_2), \varphi_o(u_2)\}$ are disjoint,
then color $v_1$ with $\varphi(u_2)$ and $v_4$ with $\varphi(u_1)$.
Now, we can greedily color $v_2$ and $v_3$ to obtain an odd 4-coloring of $G$.

(iv)
Suppose to the contrary that $G$ has a path $u_1v_1v_2u_2$ where $v_1v_2$ is a $2$-thread and $u_1$ is a vertex of odd degree. 
Let $S=\{v_1, v_2\}$, and obtain an odd 4-coloring $\varphi$ of $G-S$.
Note that $u_1\neq u_2$ by (ii).
Now, we can obtain an odd 4-coloring of $G$ by coloring $v_2$ with a color not in $\{\varphi(u_2), \varphi_o(u_2), \varphi(u_1)\}$ and $v_1$ with a color not in $\{\varphi(v_2), \varphi(u_2), \varphi(u_1)\}$.

(v)
Suppose to the contrary that $G$ has a 3-vertex $v$ that is adjacent with only $2$-neighbors $u_1, u_2, u_3$.
Let $u'_i$ be the neighbor of $u_i$ that is not $v$ for each $i\in\{1,2,3\}$. 
Let $S=\{v, u_1, u_2, u_3\}$, and obtain an odd 4-coloring $\varphi$ of $G-S$. 
Now, we can obtain an odd 4-coloring of $G$ by coloring $v$ with a color not in $\{\varphi(u'_1), \varphi(u'_2), \varphi(u'_3)\}$, and then coloring $u_i$ with a color not in $\{\varphi(u'_i), \varphi_o(u'_i), \varphi(v)\}$ for each $i$ one by one. 
\end{proof}

\begin{claim}\label{clm:mad_threads}
Let $v$ be a $4^+$-vertex of $G$.  
The following are true:
\begin{enumerate}[(i)]
    \item If $v$ has only $2$-neighbors, then $v$ is adjacent to at most $(d(v)-2)$ $2$-threads. 
    \item If $v$ is adjacent to a $3$-thread, then $v$ is adjacent to at most $d(v)-4+\min\{n_{3^+}(v),1\}$ other $2$-threads.
\end{enumerate}
\end{claim}
\begin{proof}
Let $v$ be a $d$-vertex, and let $u_1, \ldots, u_{d}$ be the neighbors of $v$.
If $u_i$ is a $2$-vertex, then let $x_i$ be the neighbor of $u_i$ that is not $v$, and for a $2$-thread $u_ix_i$, let $y_i$ be the neighbor of $x_i$ that is not $u_i$. Note that $v$, $u_i$'s, and $x_i$'s  are all distinct when $u_i$'s and $x_i$'s are $2$-vertices, since $G$ has no induced $5$-cycles, and by Claim~\ref{clm:mad_rc}~(ii).

(i) Suppose to the contrary that $u_ix_i$ is a $2$-thread for all $i\in\{2, \ldots, d\}$. We obtain an odd 4-coloring $\varphi$ of $G-S$, where $S=\{v, u_1, \ldots, u_{d}, x_2, \ldots, x_{d}\}$.
Since $G$ has no induced $5$-cycles, $x_1\neq y_i$ for all $i$, so $\varphi_o(x_1)$ is defined.
Use the color $\varphi_o(x_1)$ on $v$.
For $i\in\{2, \ldots, d\}$, use a color not in $\{\varphi(y_i), \varphi_o(y_i), \varphi(v)\}$ on $x_i$ one by one, then use a color not in $\{\varphi(x_i), \varphi(y_i), \varphi(v)\}$ on $u_i$
 for each $i\in\{2, \ldots, d\}$. 
At this point, all neighbors of $v$ except $u_1$ are colored, so $\varphi_o(v)$ is defined. 
Now, to complete an odd 4-coloring of $G$, use a color not in $\{\varphi(v), \varphi_o(v), \varphi(x_1)\}$ on $u_1$.
This is a contradiction.

(ii) Let $u_1x_1y_1$ be a $3$-thread, and let $z_1$ be the neighbor of $y_1$ that is not $x_1$.   
Let $S= \{v,y_1\}\cup N_2(v) \cup\{x_i\mid u_ix_i\text{ is a $2$-thread}\}$, and obtain an odd $4$-coloring $\varphi$ of $G-S$.  

Suppose to the contrary that $v$ is adjacent to $d-3+\min\{n_{3^+}(v), 1\}$ $2$-threads besides $u_1x_1$.
Then either $u_ix_i$ is  a $2$-thread for every $i\in\{2,\ldots, d-1\}$ and $u_d$ is the only $3^+$-neighbor of $v$, or $u_ix_i$ is a $2$-thread for every $i\in\{2,\ldots, d-2\}$ and $v$ has only $2$-neighbors. 
In the former case, let $X=\{\vp(u_d), \uc(u_d)\}$, and in the latter case, let $X=\{\vp(x_{d-1}), \vp(x_d)\}$.
Let $\alpha=\uc(z_1)$ if $\uc(z_1)$ is defined, and otherwise let $\alpha=\vp(z_1)$. 
In the latter case, use a color not in $X\cup\{\vp(z_1)\}$ on $y_1$, so the odd color of $z_1$ becomes $\vp(y_1)$. 

Now, use a color not in 
$X\cup\{\alpha\}$ on $v$, then color the  vertices in $S\setminus\{v,u_1,x_1,y_1\}$ as follows.
For each $2$-thread $u_ix_i$ where $i\ge 2$, use a color not in $\{\vp(y_i), \uc(y_i),\vp(v)\}$ on $x_i$, and
for each $2$-vertex $u_i$ where $i\ge 2$,  use a color not in $\{\vp(v), \vp(x_i), \uc(x_i)\}$ on $u_i$. 
Thus, every neighbor of $v$ is colored except $u_1$, so $\uc(v)$ is defined since $v$ has even degree by Claim~\ref{clm:mad_rc}~(iv).
Note that $\uc(y_i)$ will eventually get defined during this coloring process for $i\geq 2$.

Finally, color the uncolored vertices among $u_1,x_1,y_1$  as follows. 
If $y_1$ is uncolored, then use a color not in $\{\vp(v), \uc(v), \uc(z_1)\}$ on $u_1$, then use a color not in $\{\vp(z_1), \uc(z_1), \vp(u_1)\}$ on $y_1$.
Note that $\uc(z_1)$ is not used on $v, u_1, y_1, z_1$.
Thus, use $\uc(z_1)$ on $x_1$ to extend $\vp$ to all of $G$, which is a contradiction.
If $y_1$ is colored, then $\vp(v)=\vp(y_1)$. 
Use a color not in $\{\vp(v), \uc(v)\}$ on $u_1$, then use a color not in $\{\vp(u_1),\vp(v), \vp(z_1)\}$ on $x_1$ to extend $\vp$ to all of $G$, which is a contradiction.
\end{proof}

If a $2$-vertex $u$ is on a thread adjacent to a $3^+$-vertex $v$, then $u$ is {\it close} to $v$, and $v$ is a {\it sponsor} of $u$. 

\begin{claim}\label{clm:max_numclose}
Let $v$ be a $4^+$-vertex of $G$.
If $v$ has odd degree, then it has at most $d(v)$ close $2$-vertices. 
If $v$ has even degree, then it has at most $3d(v)-5$  close $2$-vertices.
\end{claim}
\begin{proof}
If $v$ has odd degree, then it is not adjacent to a $2$-thread by Claim~\ref{clm:mad_rc}~(iv), so it has at most $d(v)$ close $2$-vertices. 

Now assume $v$ has even degree.
If $v$ has a $3^+$-neighbor, then $v$ has at most $\max\{3(d(v)-2)+1, 2(d(v)-1)\} = 3d(v)-5$ close $2$-vertices by Claim~\ref{clm:mad_rc}~(iii) and Claim~\ref{clm:mad_threads} (ii).
If $v$ has only $2$-neighbors, then $v$ is adjacent to at most $\max\{3(d(v)-3)+3, 2(d(v)-2)+2\}=3d(v)-6$ close $2$-vertices by Claim~\ref{clm:mad_rc}~(iii) and Claim~\ref{clm:mad_threads}.
Hence, $v$ has at most $3d(v)-5$ close $2$-vertices.
\end{proof}

Let $\mu(v)$ be the initial charge of a vertex $v$, and let $\mu(v)=d(v)$. 
Let $\mu^*(v)$ be the final charge of $v$ after the following discharging rule:

\begin{enumerate}[(R1)]
    \item\label{mad_rule:1} 
    Each $3^+$-vertex sends charge $\frac{2}{9}$ to each of its close $2$-vertices.  
\end{enumerate}

Note that the total charge is preserved during the discharging process. 
Consider a vertex $v$.
We will show that $\mu^*(v) \geq \frac{22}{9}$.
Note that $G$ has no $1$-vertex by Claim~\ref{clm:mad_rc}~(i).
If $v$ is a $2$-vertex, then it receives charge $\frac{2}{9}$ from each of its two sponsors by the rule.
Hence, $\mu^*(v)=2+ \frac{2}{9} \cdot 2=\frac{22}{9}$.

If $v$ is a $3$-vertex, then at most two of its neighbors are $2$-vertices by Claim~\ref{clm:mad_rc}~(v).
Moreover, $v$ cannot be adjacent to a $2$-thread by Claim~\ref{clm:mad_rc}~(iv).
Thus, $v$ has at most two close $2$-vertices, so it sends charge $\frac{2}{9}$ at most twice by the rule.
Hence, $\mu^*(v) \ge 3- \frac{2}{9}\cdot2=\frac{23}{9}> \frac{22}{9}$.

If $v$ is a $4^+$-vertex of even degree, then it has at most $3d(v)-5$ close neighbors by Claim~\ref{clm:max_numclose}.
So it sends charge $\frac{2}{9}$ at most $3d(v)-5$ times by the rule.
Hence, $\mu^*(v) \ge d(v)-\frac{2}{9}\cdot(3d(v)-5)\geq\frac{22}{9}$ since $d(v)\geq 4$.

If $v$ is a $5^+$-vertex of odd degree, then it has at most $d(v)$ close neighbors by Claim~\ref{clm:max_numclose}.
So it sends charge $\frac{2}{9}$ at most $d(v)$ times by the rule. 
Hence, $\mu^*(v) \ge d(v)-\frac{2}{9}\cdot d(v)> \frac{22}{9}$ since $d(v)\geq 5$.

The initial charge sum is less than $\frac{22}{9}|V(G)|$, yet the final charge sum is at least $\frac{22}{9}|V(G)|$ by the above argument, so this is a contradiction. 

\section{Future directions}
\label{sec:future}

We highlight the conjecture of Petru\v sevski and \v Skrekovski~\cite{petrusevski2021colorings} regarding the odd chromatic number of planar graphs. 

\begin{conjecture}[\cite{petrusevski2021colorings}]
Every planar graph is odd 5-colorable. 
\end{conjecture}

A popular research direction is to consider planar graphs with girth restrictions.

\begin{question}
For each positive integer $c$, determine the minimum girth $g_c$ such that a planar graph with girth at least $g_c$ is odd $c$-colorable. 
\end{question}

So far, we know $6\leq g_4\leq 11$, $g_5\leq 7$, $g_7\leq g_6\leq 5$, and $g_8\leq 3$. 
Note that $g_3$ does not exist since cycles of length $\ell$ are not odd $3$-colorable when $\ell$ is not divisible by $3$, and it was conjectured that $g_4=6$ in \cite{caro2022pcfremarks}. 

Regarding odd $c$-colorings of graphs with bounded maximum average degree, the case of $c=4$ seems like an outlier. 
It would be interesting to see if one could push the threshold on the maximum average degree higher in Theorem~\ref{thm:mad22/9}.
This might be an approach to improve the upper bound on $g_4$. 
Perhaps a $5$-cycle is the only obstacle for planar graphs to be
odd $4$-colorable. 
It would be brave of us to conjecture that all planar graphs without 5-cycles are odd $4$-colorable.

\section*{Acknowledgements}
The authors would like to thank Prof. Tao Wang (Henan University) for providing comments of the previous version of this article.
Eun-Kyung Cho was supported by Basic Science Research Program through the National Research Foundation of Korea (NRF) funded by the Ministry of Education (NRF-2020R1I1A1A01058587).
Ilkyoo Choi was supported by the Basic Science Research Program through the National Research Foundation of Korea (NRF) funded by the Ministry of Education (NRF-2018R1D1A1B07043049), and also by the Hankuk University of Foreign Studies Research Fund.
Boram Park and Hyemin Kwon were supported by Basic Science Research Program through the National Research Foundation of Korea (NRF) funded by the Ministry of Science, ICT and Future Planning (NRF-2022R1F1A1069500).

\end{document}